\documentclass{article}
 \usepackage{amsmath}

\newtheorem{thm}{Theorem}[section]
\newtheorem{cor}[thm]{Corollary}
\newtheorem{prop}[thm]{Proposition}

\newtheorem{lem}[thm]{Lemma}

\newtheorem{rem}[thm]{Remark}

\newcommand{\be}{\begin{equation}}
\newcommand{\ee}{\end{equation}}
\newcommand{\ben}{\begin{enumerate}}
\newcommand{\een}{\end{enumerate}}
\newcommand{\beq}{\begin{eqnarray}}
\newcommand{\eeq}{\end{eqnarray}}
\newcommand{\beqn}{\begin{eqnarray*}}
\newcommand{\eeqn}{\end{eqnarray*}}

\newcommand{\pa}{\partial}

\newcommand{\qed}{\hspace*{\fill}Q.E.D.}  

\begin{document}
\title{On $m$-Kropina Finsler Metrics of Scalar Flag Curvature }
\author{Guojun Yang\footnote{Supported by the
National Natural Science Foundation of China (11471226) }}
\date{}
\maketitle
\begin{abstract}

 In this
paper, we consider a special class of singular Finsler metrics:
$m$-Kropina metrics which are defined by a Riemannian metric  and
a $1$-form. We show that an $m$-Kropina metric ($m\ne -1$) of
scalar flag curvature must be locally Minkowskian in dimension
$n\ge 3$. We  characterize by some PDEs a Kropina metric ($m=-1$)
which is respectively of scalar flag curvature and locally
projectively flat in dimension $n\ge 3$, and obtain some
principles and approaches of constructing non-trivial examples of
Kropina metrics of scalar flag curvature.


\end{abstract}

\section{Introduction}
The flag curvature in Finsler geometry is a natural extension of
the sectional curvature in Riemannian geometry, and  every
two-dimensional Finsler metric is of scalar flag curvature. It is
the Hilbert's Fourth Problem to study and classify projectively
flat metrics. The Beltrami Theorem  states that a Riemannian
metric is locally projectively flat if and only if it is of
constant sectional curvature. It is known that every locally
projectively flat Finsler metric is of scalar flag curvature.
However, the converse is not true.  There are regular or singular
Finsler metrics of constant flag curvature which are not locally
projectively flat (\cite{BRS} \cite{YO}). Therefore, it is an
interesting point to study and classify Finsler metrics of scalar
flag curvature. This problem is far from being solved for general
Finsler metrics. Thus we shall investigate some special classes of
Finsler metrics. Recent studies on this problem are concentrated
on Randers metrics, square metrics and some other special
$(\alpha,\beta)$-metrics.

A Randers metric is defined by $ F= \alpha+\beta$,
 where $\alpha$ is a
Riemannian metric  and
 $\beta$ is a 1-form with $b=\|\beta\|_{\alpha}<1$.
After many mathematicians' efforts, Bao-Robles-Shen finally
classify Randers metrics of constant flag curvature  by using the
navigation method (\cite{BRS}). Further, Shen-Yildirim
characterize Randers metrics of scalar flag curvature and classify
Randers metrics of weakly isotropic flag curvature (\cite{SYi}).
There are Randers metrics
 of scalar flag curvature which are neither of weakly isotropic flag curvature nor  locally
 projectively flat (\cite{CZ} \cite{SX}).   So far,  the problem of classifying Randers metrics of scalar flag
curvature still remains open.

 A square metric is written as
$ F =(\alpha+\beta)^2/\alpha$, where $\alpha$ is a Riemannian
metric and $\beta$ is a $1$-form with $b=\|\beta\|_{\alpha}<1$. In
\cite{SY}, Shen-Yildirim determine the local structure of  locally
projectively flat square metrics of constant flag curvature. Zhou
shows  that a square metric of constant flag curvature is locally
projectively flat (\cite{Z}). Later on, we  prove that a square
metric in dimension $n\ge 3$ is of scalar flag curvature iff. it
is locally projectively flat (\cite{SYa}).

In \cite{Y2}, we consider an $(\alpha,\beta)$-metric
 $F=\alpha\phi(\beta/\alpha)$ with
 $\phi(s)$ satisfying
 $$
    \big\{1+(k_1+k_3)s^2+k_2s^4\big\}\phi''(s)=(k_1+k_2s^2)\big\{\phi(s)-s\phi'(s)\big\},
 $$
where $k_1,k_2,k_3$ are constant with $k_2\ne k_1k_3$.
 We prove  that if $\beta$ is closed and the dimension $n\ge 3$, then $F$ is of scalar flag
curvature if and only if $F$ is locally projectively flat,  and
for a special case given by $\phi(s)=1+a_1s+\epsilon s^2$ with
$a_1$ and $\epsilon \ne0$ being constant, we show that $F$ is of
scalar flag curvature if and only if $F$ is locally projectively
flat.

 The Finsler metrics mentioned above  are
 regular. It seems hard to characterize a general regular
 $(\alpha,\beta)$-metric of scalar flag curvature in dimension $n\ge3$. On the other hand, singular Finsler
 metrics, such as Kropina metrics and
$m$-Kropina metrics, have a lot of applications in the real word.
In this paper, we will study $m$-Kropina metrics of scalar flag
curvature. An $m$-Kropina metric has the following form
 $$F=\alpha^{1-m}\beta^m, \ \ \ m\ne 0,1.$$
When $m=-1$, $F$ is called a Kropina metric (\cite{Kr}). There
have been  a few research papers on Kropina metrics (\cite{RR}
\cite{SYa1} \cite{YN}, \cite{Y1}--\cite{YO}). $m$-Kropina metrics
naturally appear in characterizing a class of singular
$(\alpha,\beta)$-metrics which are locally projectively flat
(\cite{Y1} \cite{Y3}) and locally projectively flat with constant
flag curvature (\cite{Y4}). Note that due to the deformation
(\ref{y11}) below for an $m$-Kropina metric, we can always assume
 $b=||\beta||_{\alpha}=1$ without loss of generality.

\begin{thm}\label{th1}
 Let $F=\alpha^{1-m}\beta^m$ be an $n(\ge 3)$-dimensional   $m$-Kropina
metric ($m\ne -1$) with $||\beta||_{\alpha}=1$. Then $F$
    is  of scalar flag curvature iff.  $F$ is locally
    Minkowskian, or more precisely, $F$ is flat-parallel ($\alpha$
    is locally flat and $\beta$ is parallel with respect to
    $\alpha$).
\end{thm}

In \cite{SYa1}, we show that an  $n(\ge 2)$-dimensional
$m$-Kropina metric ($m\ne -1$) of constant flag curvature is
locally   Minkowskian. In \cite{Y3}, we prove that  an  $n(\ge
3)$-dimensional locally projectively flat $m$-Kropina metric
($m\ne -1$) is locally
    Minkowskian.  Therefore, Theorem \ref{th1}
generalizes the corresponding results in \cite{SYa1} \cite{Y3}.
Besides,  we indicate that a two-dimensional Douglas $m$-Kropina
metric ($m\ne -1$) is  locally   Minkowskian (\cite{Y1}).

The case $m=-1$ will be much more complicated. In Section
\ref{sec4} below, we give respective  characterizations by some
PDEs for a Kropina metric to be of scalar flag curvature and
locally projectively flat in dimension $n\ge 3$ (see Theorem
\ref{th41} and Theorem \ref{prop40} below). In Section \ref{sec5},
we use Theorem \ref{th41} to prove the known local classification
for a Kropina metric of constant flag curvature (see Corollary
\ref{cor61}). However, it is difficult to determine the local
structure of a Kropina metric of scalar flag curvature, even if it
is  locally projectively flat (cf. \cite{Y1} \cite{Y3}). Here we
will show some methods (including using Corollary \ref{prop41}
below) of constructing non-trivial Kropina metrics of scalar flag
curvature.

Kropina metrics are related to Randers metrics to some extent.
Every Kropina metric is the limit of a family of Randers metrics
$F=\alpha+\beta$  as the norm $b=||\beta||_{\alpha}\rightarrow
1^{-}$ (see Remark \ref{rem72} below). Further, we have the
following result.

\begin{thm}\label{th2}
Let $F=\alpha+\beta$ be a Randers metric and $(h,\rho)$ be the
navigation data of $F$. Suppose
$\widetilde{\alpha}:=\lim_{b\rightarrow
 1^{-}}h$ is a Riemann metric and $\widetilde{\beta}:=\lim_{b\rightarrow
 1^{-}}\rho$ is a non-zero 1-form. Let $\widetilde{F}=\widetilde{\alpha}^2/\widetilde{\beta}$
 be the Kropina metric derived from $F$. Then we have
 \ben
 \item[{\rm (i)}] $\widetilde{F}=\lim_{b\rightarrow
 1^{-}}2F$, and $||\widetilde{\beta}||_{\widetilde{\alpha}}=1$.

\item[{\rm (ii)}] If $F$ is of scalar flag curvature (resp.
locally projectively flat, or Douglassian), then $\widetilde{F}$
is also of scalar flag curvature (resp. locally projectively flat,
or Douglassian). If $F$ is of weakly isotropic flag curvature,
then $\widetilde{F}$ is of constant flag curvature.
 \een

\end{thm}

For a given Randers metric $F=\alpha+\beta$ in Theorem \ref{th2},
to obtain the Kropina metric $\widetilde{F}$, we only need to
require that $\lim_{b\rightarrow
 1^{-}}(1-b^2)(\alpha^2-\beta^2)$ is a Riemann metric and $\widetilde{\beta}:=\lim_{b\rightarrow
 1^{-}}(1-b^2)\beta$ is a non-zero 1-form, since
 $h=\sqrt{(1-b^2)(\alpha^2-\beta^2)}$ and $\rho=-(1-b^2)\beta$.
 Equivalently, we can obtain $\widetilde{\alpha}$ and
 $\widetilde{\beta}$ by letting $h^{ij}\rho_i\rho_j=1$, where we put
 $h=\sqrt{h_{ij}y^iy^j}$ and $\rho=\rho_iy^i$.
By Theorem \ref{th2},  to construct non-trivial Kropina metrics of
scalar flag curvature in dimension $n\ge 3$, we can use the known
examples of Randers metrics of scalar flag curvature (see
\cite{CZ} \cite{SX}).

Next we give another principle of constructing Kropina metrics of
scalar flag curvature.

\begin{thm}\label{th3}
Let $F=\alpha^2/\beta$ be a Kropina metric with
$||\beta||_{\alpha}=1$ and define $\widetilde{F}=F+\eta$, where
$\eta$ is a closed 1-form with $||\eta||_{\alpha}$ sufficiently
small.
 \ben
 \item[{\rm (i)}] If $F$ is of scalar flag curvature, then
$\widetilde{F}$ is also a Kropina metric of scalar flag curvature.

\item[{\rm (ii)}] Let $F$ be of constant flag curvature. Then
$\widetilde{F}$ is locally projectively flat if and only if $F$ is
flat-parallel, or equivalently, $\widetilde{F}$ can be locally
written in the form
 \be\label{y5}
\widetilde{F}=\frac{|y|}{y^1}+\eta.
 \ee
 \een
\end{thm}

 Theorem \ref{th3} (ii) easily follows from a result in
\cite{Y4}, since therein we prove that
 a locally projectively flat Kropina metric with constant flag curvature is flat-parallel.
 By Theorem \ref{th3} (ii),  we can easily obtain a family of
Kropina metrics which are of scalar flag curvature but are
 neither locally projectively flat nor of constant flag
curvature in general.
 Take $\eta=\langle x,y\rangle$ with $x$ close to origin,
and then $\widetilde{F}$ in (\ref{y5}) is a projectively flat
Kropina metric with the flag curvature given by
 $$ K=\frac{3}{4}\frac{|y|^4(y^1)^4}{(\eta y^1+|y|^2)^4}.$$

 Additionally, using Corollary \ref{prop41} below and a warped product method, we obtain a family of
 Kropina metrics  which are locally projectively flat
 (see Proposition \ref{prop51} below).

\section{Preliminaries}

For a Finsler metric $F$, the Riemann curvature $R_y=R^i_{\
k}(y)\frac{\pa}{\pa x^i}\otimes dx^k$ is defined by
 \be\label{y7}
 R^i_{\ k}:=2\frac{\pa G^i}{\pa x^k}-y^j\frac{\pa^2G^i}{\pa x^j\pa
 y^k}+2G^j\frac{\pa^2G^i}{\pa y^j\pa y^k}-\frac{\pa G^i}{\pa y^j}\frac{\pa G^j}{\pa
 y^k},
 \ee
 where the spray coefficients $G^i$ are given by
\beq \label{G1}
 G^i:=\frac{1}{4}g^{il}\big \{[F^2]_{x^ky^l}y^k-[F^2]_{x^l}\big \}.
 \eeq
The Ricci curvature $Ric$ is the trace of the Riemann curvature,
that is, $ Ric:=R^m_{\ m}$. A Finsler metric is said to be of
scalar flag curvature if there is a function $ K= K(x,y)$ such
that
 \be\label{yR}
  R^i_{\ k}=KF^2(\delta^i_k-F^{-2}y^iy_k), \ \ y_k:=(F^2/2)_{y^iy^k}y^i.
  \ee
If $K$ is a constant, $F$ is said to be of constant flag
curvature. A Finsler metric $F$ is said to be projectively flat in
$U$, if
 there is a local coordinate system $(U,x^i)$ such that
  $G^i=Py^i$,
   where
$P=P(x,y)$ is called the projective factor satisfying $P(x,\lambda
y)=\lambda P(x,y)$ for $\lambda>0$.

The Weyl curvature $W^i_{\ k}$ and the Douglas curvature $D^{\
i}_{h \ jk}$ are two important projectively invariant tensors and
they are defined respectively by
 \beq\label{y8}
 W^i_{\ k}:\hspace{-0.6cm}&&=R^i_{\ k}-\frac{R^m_{\
 m}}{n-1}\delta^i_k-\frac{1}{n+1}\frac{\pa }{\pa
 y^m}\big(R^m_{\ k}-\frac{R^h_{\
 h}}{n-1}\delta^m_k\big)y^i,\\
D^{\ i}_{h \ jk}:\hspace{-0.6cm}&&=\frac{\pa^3}{\pa y^h\pa y^j\pa
 y^k}\big(G^i-\frac{1}{n+1}\frac{\pa G^m}{\pa
 y^m}y^i\big).\nonumber
 \eeq
 In two-dimensional case, there is a projectively invariant
 tensor $W^o$ called Berwald-Weyl tensor.
 A Finsler metric is called a Douglas metric if $D^{\
i}_{h \ jk}=0$. A Finsler metric is of scalar flag curvature if
and only if $W^i_{\ k}=0$. An $n$-dimensional Finsler metric is
locally projectively flat if and only if: $W^i_{\ k}=0$ and $D^{\
i}_{h \ jk}=0$ for $n\ge 3$, and $W^o=0$ and $D^{\ i}_{h \ jk}=0$
for $n=2$ (\cite{Ma1}).

An $(\alpha,\beta)$-metric $F$ is a Finsler metric  defined by a
Riemannian metric
 $\alpha=\sqrt{a_{ij}(x)y^iy^j}$ and a $1$-form $\beta=b_i(x)y^i$ on a manifold
 $M$, which is expressed in the following form:
 $$F=\alpha \phi(s),\ \ s=\beta/\alpha,$$
where $\phi(s)$ is a suitable function. If we take $\phi(s)=1+s$,
then we get the well-known Randers metric $F=\alpha+\beta$. In
applications, there are a lot of singular Finsler metrics. In this
paper, we will discuss a class of singular $(\alpha,\beta)$
Finsler metrics---$m$-Kropina metrics.

 An $m$-Kropina metric
 is in the form  $F=\alpha^{1-m}\beta^m$, where
$m\ne 0,1$ is real. In particular, it is called a Kropina metric
when $m=-1$. For an $m$-Kropina metric $F=\alpha^{1-m}\beta^m$, we
introduce  a special deformation on $\alpha$ and $\beta$.  Define
a new pair $(\widetilde{\alpha},\widetilde{\beta})$ by
 \be\label{y11}
\widetilde{\alpha}:=b^m\alpha,\ \ \
\widetilde{\beta}:=b^{m-1}\beta,
 \ee
which appears first in \cite{SYa1}. It is interesting that under
the deformation (\ref{y11}), the $m$-Kropina metric
$F=\alpha^{1-m}\beta^m$ keeps formally unchanged, that is,
 \be\label{y011}
 F=\widetilde{\alpha}^{1-m}\widetilde{\beta}^m,\ \ \ \
 (||\widetilde{\beta}||_{\widetilde{\alpha}}=1).
 \ee
It has been shown  that the deformation (\ref{y11}) plays an
important role  on the study of $m$-Kropina metrics (\cite{SYa1}
\cite{Y1}--\cite{Y4}). Due to (\ref{y011}), we can always assume
$||\beta||_{\alpha}=1$ for an $m$-Kropina metric
$F=\alpha^{1-m}\beta^m$ without loss of generality.

 For a Riemann metric  $\alpha =\sqrt{a_{ij}y^iy^j}$ and a
$1$-form $\beta = b_i y^i $, define
 $$r_{ij}:=\frac{1}{2}(b_{i|j}+b_{j|i}),\ \ s_{ij}:=\frac{1}{2}(b_{i|j}-b_{j|i}),\ \
 r^i_{\ j}:=a^{ik}r_{kj},\ \  s^i_{\ j}:=a^{ik}s_{kj},$$
 $$p_{ij}:=r_{ik}r^k_{\ j},\ \  q_{ij}:=r_{ik}s^k_{\ j}, \ \ t_{ij}:=s_{ik}s^k_{\ j},\ \
 r_j:=b^ir_{ij},\ \  s_j:=b^is_{ij},$$
 $$
  p_j:=b^ip_{ij},\ \  q_j:=b^iq_{ij}, \ \ r_j:=b^ir_{ij},\ \ t_j:=b^ir_{ij},\ \ r:=b^ir_i,
 $$
 where $b^i$ is defined by $b^i:=a^{ij}b_j$, $(a^{ij})$ is the inverse of
 $(a_{ij})$, and $\nabla \beta = b_{i|j} y^i dx^j$  denotes the covariant
derivatives of $\beta$ with respect to $\alpha$.  We use $a_{ij}$
to raise or lower the indices of a tensor. For a  tensor $T_{ij}$
as an example, define $T_{i0}:=T_{ij}y^j$ and
$T_{00}:=T_{ij}y^iy^j$, etc.

\begin{lem}\label{lem21}
  Here we list some identities as follows:
  \beq
  &&q_{ik}+q_{ki}=r_{i|k}+r_{k|i}-2p_{ik}-b^m(r_{mi|k}+r_{mk|i}),\  \ \ \ \ \ q_{ik}=t_{ik}+s_{k|i}-b^ms_{mk|i}, \label{j8}\\
  &&s_{ij|k}=r_{ik|j}-r_{jk|i}-b_l\bar{R}^{\ l}_{k\ ij},\ \ b^mb^v(r_{mv|k}-r_{mk|v})=t_k-q_k+b^ms_{k|m}, \label{j08}\\
  &&b^mq_{km}=-r_{km}s^m=b^ms_{m|k}+t_k, \ \ \ \ \ \  q^k_{\ k}=0,\label{j008}
 \eeq
  where $\bar{R}$ denotes the Riemann curvature tensor of
  $\alpha$. If $||\beta||_{\alpha}=constant$,
we have
 \beq
  &&  r_k+s_k=0,\ \ \  b^mq_m=s_ms^m,\ \  \ b^lb^ks_{k|l}=2s^ls_l,\ \ \ b^ib^jb^kr_{ij|k}=-4s_ls^l,\label{j9}\\
  &&    b^ms_{i|m}=2s^mr_{im}-b^kb^lr_{ki|l}=-2b^mq_{im}-b^kb^lr_{ki|l}.\label{j10}
  \eeq
\end{lem}

For an $m$-Kropina metric $F=\alpha^{1-m}\beta^m$, by (\ref{G1})
we get
 \be\label{y12}
 G^i=G^i_{\alpha}-\frac{m}{(m-1)s}\alpha
 s^i_0+\frac{m}{2(m-1)}\frac{(m-1)sr_{00}+2m\alpha s_0}{
 s\big[mb^2-(m+1)s^2\big]}(b^i-2\alpha^{-1}sy^i).
 \ee
Then by (\ref{y8}) and (\ref{y12}), we can get the expressions of
the Weyl curvature tensor $W^i_{\ k}$ for an $m$-Kropina metric
$F=\alpha^{1-m}\beta^m$. We have given a Maple program in
\cite{SYa} to compute
 the Weyl curvature for any $(\alpha,\beta)$-metric. In this
 paper, we will write out the whole expression of the Weyl
 curvature for a Kropina metric ($m=-1$); while for $m\ne-1$, we will not write out the
 expression since it is very long, but some key terms will be given, similarly
like what we have done in studying square metrics in \cite{SYa}.

For readers to verify the expression of $W^i_{\ k}$ for  an
$m$-Kropina metric $F=\alpha^{1-m}\beta^m$, we give the expression
of a leading term. We see that $W^i_{\ k}\times
(n^2-1)(m-1)^2\beta^3\big[mb^2\alpha^2-(m+1)\beta^2\big]^5=0$ can
be written as
 \be\label{weyl}
(n+1)m^7b^8A_{14}\alpha^{14}+A_{12}\alpha^{12}+A_{10}\alpha^{10}+\cdots+A_2\alpha^2+A_0=0,
 \ee
where $A_0,A_2,\cdots, A_{14}$ are polynomials in $(y^i)$, and
$A_{14}$ is given by
 \beqn
A_{14}\hspace{-0.6cm}&&=-(n-1)\big[(b^2t^i_{\
0}+s_0s^i)b_k-(t_0b_k+\beta t_k)b_i+\beta(b^2t^i_{\
k}+s^is_k)\big]\\
&&\ \ \   +(2s_js^j+b^2t^j_{\ j})(\beta\delta^i_k+y^ib_k).
 \eeqn
When $m=-1$, eliminating the factor $-b^6\alpha^{10}$ from
(\ref{weyl}) we obtain
 \be\label{K10}
 (n+1)b^2B_4\alpha^4+2(n+1)\beta B_2\alpha^2+4\beta^2B_0=0,
 \ee
where $B_4=A_{14}$, and $B_2,B_0$ are given by (denote by
$\bar{W}^i_{\ k}$ the Weyl curvature of $\alpha$)
 \beqn
B_2\hspace{-0.6cm}&&=(b^4s^j_{\
0|j}+b^2q_0-b^2b^js_{0|j}-b^2r^j_{\
j}s_0-b^2b^jq_{0j}+rs_0)(2\beta\delta^i_k+y^ib_k)\\
&&-y^i\big[b^2(2s_js^j+b^2t^j_{\ j})y_k+\beta(r-b^2r^j_{\
j})s_k+b^2\beta(q_k+b^2s^j_{\ k|j}-b^jq_{kj}-b^js_{k|j})\big]\\
&&+(n-1)\big[b^2(b^2t^i_{\
0}+s^is_0-t_0b^i)y_k+(b^2s_{0|0}-b^2q_{00}-r_0s_0-s_0^2)b^ib_k+\beta(r_0-s_0)b^is_k\\
&&+\beta
b^i(b^2q_{0k}-2b^2q_{k0}+2b^2s_{0|k}-b^2s_{k|0}-2s_0r_k)-(b^2s^i_{\
0|0}+r_{00}s^i-r^i_{\ 0}s_0)b_k\\
&&-b^2\beta r^i_{\ 0}s_k-b^2\beta(r_{k0}s^i-2r^i_{\
k}s_0+2b^2s^i_{\ 0|k}-b^2s^i_{\ k|0})\big],\\
   B_0\hspace{-0.6cm}&&=(n+1)\beta(b^2r_{0|0}-b^2r^j_{\
   j}r_{00}+b^2s_{0|0}-2b^2q_{00}-b^2b^jr_{00|j}-r_0^2-s_0^2+rr_{00}-2r_0s_0)\delta^i_k\\
   &&-(n+1)(b^4s^j_{\
0|j}+b^2q_0-b^2b^js_{0|j}-b^2r^j_{\
j}s_0-b^2b^jq_{0j}+rs_0)y^iy_k\\
&&+\beta y^i\big[(n+1)(r_0+s_0)(r_k+s_k)-(n+1)(r-b^2r^j_{\
j})r_{k0}+(n-2)b^2r_{k|0}-(2n-1)b^2r_{0|k}\\
&&+(n+1)b^2(q_{k0}+q_{0k})+(n-2)b^2s_{k|0}-(2n-1)b^2s_{0|k}+(n+1)b^2b^jr_{k0|j}\big]\\
&&+(n^2-1)\big\{y_k[(r_0s_0+s_0^2+b^2q_{00}-b^2s_{0|0})b^i+b^2r_{00}s^i+b^2(b^2s^i_{\
0|0}-s_0r^i_{\ 0})]\\
&&+\beta
b^i[(r_0+s_0)r_{k0}-(r_k+s_k)r_{00}-b^2(r_{k0|0}-r_{00|k})]+b^2\beta(b^2\bar{W}^i_{\
k}+r_{00}r^i_{\ k}-r^i_{\ 0}r_{k0})\big\}.
 \eeqn

\begin{lem}
 Let $F=\alpha^2/\beta$ be an $n$-dimensional Kropina metric. Then
 $W^i_{\ k}=0$ is equivalent to (\ref{K10}), and the Ricci
 curvature $Ric$ of $F$ is given by
  \beq
 Ric\hspace{-0.5cm}&&=\bar{R}ic-\frac{1}{4b^4\alpha^2 s^2}\Big\{b^2(b^2t^l_{\
 l}+2s_ls^l)\alpha^4+4s\big[b^4s^l_{\ 0|l}-(n-1)b^2t_0+(r-b^2r^l_{\
 l})s_0\nonumber\\
 &&+b^2(q_0-b^ls_{0|l}-b^lq_{0l})\big]\alpha^3+4s^2\big[(r-b^2r^l_{\
 l})r_{00}+(n-2)s_0^2+2(2n-3)r_0s_0\nonumber\\
 &&-(n-2)b^2s_{0|0}+b^2r_{0|0}-r_0^2-2nb^2q_{00}-b^lr_{00|l}\big]\alpha^2\nonumber\\
 &&+4(n-1)s^3\big[2r_{00}(2r_0-s_0)-b^2r_{00|0}\big]\alpha-12(n-1)s^4r_{00}^2\Big\},\label{ric}
  \eeq
  where $\bar{R}ic$ denotes the Ricci curvature of $\alpha$.
\end{lem}

\section{Proof of Theorem \ref{th1}}

\begin{lem}\label{lem31}
 $\beta$ is closed $\Longleftrightarrow t_{ij}=0$ $\Longleftrightarrow t^k_{\
 k}=0$.
\end{lem}

\begin{lem}\label{lem33}
 Let  $F=\alpha^{1-m}\beta^m$ be an $m$-Kropina
metric ($m\ne -1$) of scalar flag curvature on an  $n(\ge
3)$-dimensional manifold $M$. Then $r_{00}$ satisfies
 \be\label{y14}
 r_{00}=2\tau
 \big[mb^2\alpha^2-(m+1)\beta^2\big]-\frac{2(m+1)}{(m-1)b^2}\beta
 s_0,
 \ee
 where $\tau=\tau(x)$ is a scalar function.
\end{lem}

{\it Proof :}  Since $F=\alpha^{1-m}\beta^m$ is of scalar flag
curvature, we have $W^i_{\ k}=0$. Then we have (\ref{weyl}).  Now
$\alpha^2\times (\ref{weyl})$ can be written as
 \be\label{y15}
 C^i_k\big[mb^2\alpha^2-(m+1)\beta^2\big]-24(n-2)(m+1)^3\beta^8y^i(\alpha^2
 b_k-\beta
 y_k)\big[(m-1)\beta r_{00}+2m\alpha^2 s_0\big]^2=0,
 \ee
where $C^i_k$ are polynomials in $(y^i)$. It is easy to see that
$mb^2\alpha^2-(m+1)\beta^2$ is an irreducible polynomial in
$(y^i)$ since $m\ne 0$ and $n>2$. Further, if $\alpha^2 b_k-\beta
y_k$ is divisible by $mb^2\alpha^2-(m+1)\beta^2$ for all $k$, then
there are scalar functions $\tau_k=\tau_k(x)$ such that
 $$
\alpha^2 b_k-\beta y_k=\tau_k\big[mb^2\alpha^2-(m+1)\beta^2\big].
 $$
 Contracting the above by $y^k$ we have $\tau_0=0$ and hence $\alpha^2 b_k-\beta
 y_k=0$. This is a contradiction.
Now since $n>2$ and $m\ne -1$,  it follows from (\ref{y15})
  that $(m-1)\beta r_{00}+2m\alpha^2 s_0$ is
divisible by $mb^2\alpha^2-(m+1)\beta^2$, which implies
 \be\label{y18}
(m-1)\beta r_{00}+2m\alpha^2
s_0=\theta\big[mb^2\alpha^2-(m+1)\beta^2\big],
 \ee
 where $\theta$ is a 1-form. Eq. (\ref{y18}) is equivalent to
  \be\label{y018}
 m(2s_0-b^2\theta)\alpha^2+\beta\big[(m-1)r_{00}+(m+1)\theta\beta\big]=0.
  \ee
By (\ref{y018}), there is a scalar function $\tau=\tau(x)$ such
that
 \be\label{y0018}
 2s_0-b^2\theta=-2(m-1)b^2\tau\beta.
 \ee
Now plugging (\ref{y0018}) into (\ref{y018}) immediately yields
(\ref{y14}). \qed

\begin{lem}\label{lem34}
 Let  $F=\alpha^{1-m}\beta^m$ be an $m$-Kropina
metric ($m\ne -1$) of scalar flag curvature. Then we have
 \be\label{y21}
 t^k_{\ k}=-\frac{2s_ks^k}{b^2}.
 \ee
\end{lem}

{\it Proof :} Since $F=\alpha^{1-m}\beta^m$ is of scalar flag
curvature, we have (\ref{weyl}), and further we can rewrite
(\ref{weyl})  as
 \be\label{y22}
 D^i_k\beta+m^6(n+1)b^8\alpha^{12}b_kT^i=0,
 \ee
where $D^i_k$ are polynomial in $(y^i)$ and $T^i$ are defined by
 $$T^i:=m\big[(n-1)(b^it_0-s^is_0-b^2t^i_{\ 0})+y^i(b^2t^j_{\
 j}+2s_js^j)\big]\alpha^2+2(m+1)(b^2t_{00}+s_0^2)y^i.$$
Now it follows from (\ref{y22}) that there are polynomials $f^i$
in $(y^i)$ of degree two  such that
 \be\label{y23}
 T_i-f_i\beta=0.
 \ee
Contracting (\ref{y23}) by $y^i$ we get
 \be\label{y24}
 m(2s_ks^k+b^2t^k_{\
 k})\alpha^4+\big[(2+3m-nm)(b^2t_{00}+s_0^2)+m(n-1)\beta t_0\big]\alpha^2-f_0\beta=0.
 \ee
Then by (\ref{y24}), we have $f_0=\theta \alpha^2$ for some 1-form
$\theta=\theta_i(x)y^i$. Plugging it into (\ref{y24}) gives
 \beq\label{y25}
 0&=&2m(2s_ks^k+b^2t^k_{\
 k})a_{ij}+2(2+3m-nm)(b^2t_{ij}+s_is_j)+\nonumber\\
 &&m(n-1)(b_it_j+b_jt_i)-(b_i\theta_j+b_j\theta_i).
 \eeq
Contracting (\ref{y25}) by $a^{ij}$ yields
 \be\label{y26}
 (2+3m)b^2t^k_{\ k}+2(1+2m)s_ks^k-b^k\theta_k=0.
 \ee
Further contracting (\ref{y25}) by $b^ib^j$ gives
 \be\label{y27}
 mb^2t^k_{\ k}-2s_ks^k-b^k\theta_k=0.
 \ee
Now it is easy to follow from (\ref{y26}) and (\ref{y27}) that
(\ref{y21}) holds.     \qed

\

{\it Proof of Theorem \ref{th1} :} \  Let $F=\alpha^{1-m}\beta^m$
be an $n(\ge 3)$-dimensional $m$-Kropina metric ($m\ne -1$) of
scalar flag curvature.  Then under the deformation (\ref{y11}),
$F=\widetilde{\alpha}^{1-m}\widetilde{\beta}^m$ is also an
$m$-Kropina metric of scalar flag curvature. So we obtain Lemma
\ref{lem33} and Lemma \ref{lem34} under $\widetilde{\alpha}$ and
$\widetilde{\beta}$.

Note that $\widetilde{b}^2=1$, and then by (\ref{y14}) we have
 \be\label{y28}
 \widetilde{r}_{ij}=2\widetilde{\tau}\big[m\widetilde{a}_{ij}-(m+1)\widetilde{b}_i\widetilde{b}_j\big]
 -\frac{m+1}{m-1}(\widetilde{b}_i\widetilde{s}_j+\widetilde{b}_j\widetilde{s}_i),
 \ee
We will prove $\widetilde{r}_{ij}=0$ by (\ref{y28}). This fact is
essentially proved in \cite{SYa1} \cite{Y3}. For convenience, we
give the proof here. Contracting (\ref{y28}) by $\widetilde{b}^i$
and using $||\widetilde{\beta}||_{\widetilde{\alpha}}=constant=1$
we have
 \be\label{y29}
 \widetilde{r}_j+\widetilde{s}_j=-2\widetilde{\tau}
 \widetilde{b}_j-\frac{2}{m-1}\widetilde{s}_j=0.
 \ee
Contracting (\ref{y29}) by $\widetilde{b}^j$ we get
$\widetilde{\tau}=0$ and then by (\ref{y29}) again we have
$\widetilde{s}_j=0$.  Thus by (\ref{y28}) again we have
 $$\widetilde{r}_{ij}=0.$$

Next by (\ref{y21}) we have
 \be\label{y30}
 \widetilde{t}^k_{\ k}=-2\widetilde{s}_k\widetilde{s}^k.
 \ee
Since we have proved $\widetilde{s}_k=0$, we have
$\widetilde{t}^k_{\ k}=0$ by (\ref{y30}). Thus  Lemma \ref{lem31}
implies that $\widetilde{\beta}$ is closed. Thus by this fact and
$\widetilde{r}_{ij}=0$, we obtain that $\widetilde{\beta}$ is
parallel with respect to $\widetilde{\alpha}$.   \qed

\section{Kropina metrics of scalar flag curvature}\label{sec4}

\subsection{Main results}

\begin{thm}\label{th41}
 Let $F=\alpha^2/\beta$ be an $n(\ge 2)$-dimensional Kropina
 metric with $||\beta||_{\alpha}=1$. Denote by  $\bar{R}^i_{\ k}$  the Riemann curvature tensor of
$\alpha$. Then $F$ is of scalar flag curvature if and only if the
 following hold
  \beq
  s_{ij|k}&=&\Big\{t_j-\frac{t^l_{\
  l}-(n-3)s^ls_l}{n-1}
  b_j\Big\}a_{ik}+r_{ik}s_j+q^*_{ki}b_j+s_{j|k}b_i-(i/j),\label{y032}\\
 \bar{R}^i_{\ k}&=&\frac{(n-3)s^ls_l-t^l_{\
 l}}{n-1}
 \big(\alpha^2\delta^i_k-y^iy_k\big)-B_{00}\delta^i_k-B^i_{\ k}\alpha^2
 +B_{0k}y^i+B^i_{\ 0}y_k\nonumber\\
 &&+r^i_{\
 0}r_{k0}-r_{00}r^i_{\ k},\label{y49}
  \eeq
  where the symbol $(i/j)$ above denotes the terms obtained from
the proceeding terms by the interchange of the indices $i$ and
$j$,  and $q^*_{ik}$, $\sigma_i$ and $B^i_{\ k}$ are defined by
   \beq
   q^*_{ik}:&=&\frac{1}{2}b^pb^l\big[(r_{lp|i}-r_{li|p})b_k-(i/k)\big]-\frac{1}{2}b^l(r_{lk|i}
  +r_{li|k})
  -p_{ik}-s_{i|k},\label{Q}\\
 \sigma_i:&=&2\big[(n-3)s^ls_l-t^l_{\
 l}-(n-1)\lambda\big]b_i+2(n-1)b^pb^l(r_{lp|i}-r_{li|p}),\label{j00}\\
 B_{ik}:&=&\frac{1}{2}(r_{il}r^l_{\ k}+b^lr_{lk|i})
  +\frac{b_k\sigma_i}{4(n-1)}+s_{i|k}+(i/k),\label{j01}
   \eeq
and $\lambda=\lambda(x)$ is a scalar function.
 In this case, the  flag curvature $K$ of
$F$
 is given by
 \beq\label{y050}
K&=&\lambda
s^2+\frac{s^2}{\alpha^2}\Big\{\frac{3s^2}{\alpha^2}r_{00}^2+\frac{s}{\alpha}(r_{00|0}
+6r_{00}s_0)+3q_{00}+3s_0^2-b^l(r_{l0|0}-r_{00|l})\Big\}\nonumber\\
&&+\frac{1}{4(n-1)}\big[(4s^2-1)t^l_{\
l}-2(1+2ns^2-6s^2)s^ls_l\big].
 \eeq
\end{thm}

 In \cite{Y1} \cite{Y3}, we give a  way  to
characterize locally projectively flat Kropina metrics in
dimension $n\ge 2$ by (\ref{ycw1}) and an equation on the spray
$G^i_{\alpha}$ of $\alpha$. Now using Theorem \ref{th41}, we can
obtain a different way to characterize locally projectively flat
Kropina metrics by adding a Douglasian condition (\ref{ycw1}) in
$n\ge 3$.

\begin{thm}\label{prop40}
Let $F=\alpha^2/\beta$ be an $n(\ge 3)$-dimensional Kropina
 metric with $||\beta||_{\alpha}=1$.  Then $F$ is locally projectively flat
  if and only if (\ref{y49}) and the
 following hold
  \beq
 s_{ij}=b_is_j-b_js_i.\label{ycw1}
  \eeq
  In this case,   the  flag curvature $K$ of $F$ is given by (\ref{y050}), and $\sigma_i$ in (\ref{y49}) are
 given by
 \be\label{j000}
 \sigma_i=2(n-1)\big[b^ls_{i|l}-(\lambda+s_ls^l)b_i\big].
 \ee
\end{thm}

 In a special
case, we have the following simple corollary. We will construct
some examples in Section \ref{sec6} below by  Corollary
\ref{prop41}.

\begin{cor}\label{prop41}
 Let $F=\alpha^2/\beta$ be an $n(\ge 3)$-dimensional Kropina
 metric  with $||\beta||_{\alpha}=1$. Suppose
  \be\label{y53}
 b_{i|j}=\epsilon (a_{ij}-b_ib_j),\ \ \
 \epsilon_{i}=ub_i,
  \ee
  where $u=u(x),\epsilon=\epsilon(x)$ are scalar  functions and
  $\epsilon_i:=\epsilon_{x^i}$. Then $F$ is locally projectively flat  if and only if
   \be
  \bar{R}^i_{\ k}=-\epsilon^2(\alpha^2\delta^i_k-y^iy_k)-u(\alpha^2
  b^ib_k+\beta^2\delta^i_k-\beta y^ib_k-\beta y_kb^i).\label{y54}
   \ee
  In this
   case, the  flag curvature $K$ is given by
    \be\label{y056}
K =s^6\big[\epsilon^2 (3s^2-4)-u\big].
    \ee
\end{cor}

\begin{rem}
 It is known in \cite{SYi} that  a Randers metric $F=\alpha+\beta$ in dimension
 $n\ge 2$ is of scalar flag curvature if and only if for some scalar
 $\lambda=\lambda(x)$,
   \beq
s_{ij|k}&=&\frac{1}{n-1}(a_{ik}s^m_{\ j|m}-a_{jk}s^m_{\ i|m}),\label{jj44}\\
\bar{R}^i_{\ k}&=&\lambda
 (\alpha^2\delta^i_k-y^iy_k)+\alpha^2t^i_{\
 k}+t_{00}\delta^i_k
 -t_{k0}y^i-t^i_{\ 0}y_k-3s^i_{\ 0}s_{k0}.\label{jj45}
 \eeq
 In Theorem \ref{th41}, for a Kropina metric of scalar flag curvature, we obtain the equations (\ref{y032})
and (\ref{y49}) similar to (\ref{jj44}) and (\ref{jj45}). However,
the characterization and proof for a Kropina metric are much more
complicated than that for a Randers metric.
\end{rem}

\subsection{Proof of Theorem \ref{th41}}

\begin{prop}\label{prop45}
 Let $F=\alpha^2/\beta$ be an $n(\ge 2)$-dimensional Kropina
 metric with $||\beta||_{\alpha}=1$.  Then $F$ is of scalar flag curvature if and only if  (\ref{y49}) and the
 following hold
  \beq
  t_{ij}&=&b_it_j+b_jt_i-s_is_j+\frac{1}{n-1}\Big\{(t^l_{\
  l}+2s^ls_l)a_{ij}-\big[t^l_{\
  l}-(n-3)s^ls_l\big]b_ib_j\Big\},\label{y31}\\
  s_{ij|k}&=&\Big\{t_j-\frac{t^l_{\
  l}-(n-3)s^ls_l}{n-1}
  b_j\Big\}a_{ik}+r_{ik}s_j+q_{ki}b_j+s_{j|k}b_i-(i/j),\label{y32}\\
 q_{ik}&=&\frac{1}{2}b^mb^l\big[(r_{lm|i}-r_{li|m})b_k-(i/k)\big]-\frac{1}{2}b^l(r_{lk|i}
  +r_{li|k})
  -p_{ik}-s_{i|k}.\label{y48}
  \eeq
\end{prop}

 {\it Proof :}  Assume  $F=\alpha^2/\beta$ is of scalar flag
curvature in dimension $n\ge 2$. By $W^i_{\ k}=0$, we get
(\ref{K10}). Here we put $b=||\beta||_{\alpha}=1$ and hence
$r_k+s_k=0=r$ in (\ref{K10}). First, (\ref{K10}) can be written as
 \be\label{y34}
 (\cdots)\beta-(n+1)\alpha^4b_k\big[(n-1)(t_{i0}-t_0b_i+s_0s_i)-(t^l_{\
 l}+2s^ls_l)y_i\big]=0,
 \ee
where the omitted term is a homogeneous polynomial in $(y^i)$.
Then by (\ref{y34}) we have
 \be\label{y35}
t_{ij}=t_jb_i-s_is_j-\frac{\rho_ib_j-(t^l_{\
 l}+2s^ls_l)a_{ij} }{n-1},
 \ee
where $\rho_i=\rho_i(x)$ are some scalar functions. By
(\ref{y35}), using $t_{ij}=t_{ji}$ we get
 \be\label{y36}
 \rho_i=\sigma b_i-(n-1)t_i,
 \ee
where $\sigma=\sigma(x)$ is a scalar function. Plugging
(\ref{y36}) into (\ref{y35}) and then contracting (\ref{y35}) by
$a^{ij}$, we get
 \be\label{y37}
\sigma=t^l_{\ l}-(n-3)s^ls_l.
 \ee
Now plugging (\ref{y36}) and (\ref{y37}) into (\ref{y35}) we
obtain (\ref{y31}).

By  $t_{ik}$ and $t_{i0}$ given by (\ref{y31}), we can write
$(\ref{y34})/\beta$ as
 \be\label{y38}
 (\cdots)\beta+2(+1)\alpha^2 b_kC_i=0,
 \ee
where $C_i$ is a  homogeneous polynomial of degree two in $y$ (the
expression is omitted here). It is easy to see from (\ref{y38})
that $C_i$ is divisible by $\beta$. Hence we have
$C_i=c_{i0}\beta$ for a 1-form $c_{i0}=c_{ij}y^j$, which is
equivalent to
 \beq\label{y39}
 s_{i0|0}&=&\frac{q_0+s^l_{\ 0|l}-b^lq_{0l}-b^ls_{0|l}-r^l_{\
 l}s_0}{n-1}\ y_i+\Big\{\frac{t^l_{\
 l}-(n-3)s^ls_l}{n-1}\alpha^2-q_{00}+s_{0|0}\Big\}b_i\nonumber\\
 &&-\alpha^2
 t_i+s_0r_{i0}-r_{00}s_i-\frac{c_{i0}\beta}{n-1}.
 \eeq
 Plug (\ref{y39}) into
(\ref{y38}) and then $(\ref{y38})/(2\beta)$ can be written as
 \be\label{y40}
 (\cdots)\beta+(n+1)\alpha^2D_{ik}=0,
 \ee
where $D_{ik}$ is a 1-form (the expression is omitted here). It is
easy to see from (\ref{y40}) that $D_{ik}$ is divisible by
$\beta$. Hence we have $D_{ik}=f_{ik}\beta$ for a scalar function
$f_{ik}$, which is equivalent to
 \be\label{y41}
 (n-1)s_{ik|j}-2(n-1)s_{ij|k}+\cdots=f_{ik}b_j.
 \ee
 Interchanging $j,k$
in (\ref{y41}) we have
 \be\label{y42}
(n-1)s_{ij|k}-2(n-1)s_{ik|j}+\cdots=f_{ij}b_k,
 \ee
Then $2\times (\ref{y41})+(\ref{y42})$ gives
 \beq\label{y43}
 s_{ij|k}&=&\frac{q_j+s^l_{\ j|l}-b^lq_{jl}-b^ls_{j|l}-r^l_{\ l}s_j}{n-1}\ a_{ik}+\Big\{\frac{t^l_{\
 l}-(n-3)s^ls_l}{n-1}b_i-t_i\Big\}a_{jk}\nonumber\\
 &&+\frac{2b_kc_{ij}+b_jc_{ik}-b_kf_{ij}-2b_jf_{ik}}{3(n-1)}-b_iq_{kj}+b_is_{j|k}+s_jr_{ik}-s_ir_{jk}.
 \eeq
By (\ref{y39}) and (\ref{y43}) we get
 \be\label{y44}
 f_{ij}=2c_{ij}.
 \ee
By (\ref{y44}) and $s_{ij|k}+s_{ji|k}=0$, it follows from
(\ref{y43}) that
 \beq\label{y45}
0&=&\Big\{\frac{q_j+s^l_{\ j|l}-b^lq_{jl}-b^ls_{j|l}-r^l_{\
l}s_j}{n-1}+\frac{t^l_{\ l}-(n-3)s^ls_l}{n-1}b_j-t_j\Big\}
a_{ik}\nonumber\\
&&-b_iq_{kj}+b_is_{j|k}-\frac{b_ic_{jk}}{n-1}+(i/j).
 \eeq
Contracting (\ref{y45}) by $b^ib^j$ we can first get the
expression of $b^lc_{lk}$, and then using $b^lc_{lk}$ and
contracting (\ref{y45}) by $b^j$ we can get the expression of
$c_{ik}$. Now plugging $c_{ik}$ into (\ref{y45}) yields
 \be\label{y46}
 0=\Big\{\big[\frac{b^m(b^ls_{m|l}-s^l_{\ m|l})}{n-1}-s^ls_l\big]b_j+\frac{q_j+s^l_{\ j|l}-b^lq_{jl}-b^ls_{j|l}-r^l_{\
l}s_j}{n-1}-t_j\Big\}(a_{ik}-b_ib_k)+(i/j).
 \ee
Contracting (\ref{y46}) by $a^{ik}$ we obtain
 \be\label{y47}
 s^l_{\ j|l}=b^l(q_{jl}+s_{j|l})+(n-1)t_j+r^l_{\
l}s_j-q_j+\big[(n-1)s^ls_l-b^m(b^ls_{m|l}-s^l_{\ m|l})\big]b_j.
 \ee
Finally, plugging (\ref{y44}), $c_{ij}$ and (\ref{y47}) into
(\ref{y43}) we obtain (\ref{y32}).

 By (\ref{y32}), we can
determine the expressions of the following quantities
 $$s_{ik|0},\ s^l_{\ 0|l},\ s^l_{\ k|l},\ s_{i0|k},\ s_{i0|0},\ b^ms^l_{\ m|l}.$$
Plug these  quantities into (\ref{y38}) and then (\ref{y38}) is
equivalent to ($\bar{W}_{ik}:=a_{il}\bar{W}^l_{\ k}$)
 \beq\label{y33}
 \bar{W}_{ik}&=&\frac{1}{n-1}\big\{s^l_{|l}(\alpha^2a_{ik}-y_iy_k)+(2q_{00}+b^lr_{00|l}+r^l_{\
  l}r_{00})a_{ik}-\nonumber\\
  &&(r^l_{\
  l}r_{k0}+b^lr_{k0|l}+q_{k0}+q_{0k})y_i\big\}+(s_{i|0}-q_{0i})y_k+(r_{k0|0}-r_{00|k})b_i+\nonumber\\
  &&(q_{ki}-s_{i|k})\alpha^2+r_{k0}r_{i0}-r_{00}r_{ik}.
 \eeq

\begin{lem}
(\ref{y33}) is equivalent to the following equation
 \beq
\bar{R}^i_{\ k}&=&\lambda
 (\alpha^2\delta^i_k-y^iy_k)+\big[b^l(r_{00|l}-r_{l0|0})+q_{00}-s_{0|0}\big]\delta^i_k+r^i_{\
 0}r_{k0}-r_{00}r^i_{\ k}\nonumber\\
 &&+(q^{\ i}_k-s^i_{|k})\alpha^2+\frac{1}{2}\big[b^l(r_{l0|k}+r_{lk|0}-2r_{k0|l})
 -q_{k0}-q_{0k}+s_{k|0}+s_{0|k}\big]y^i\nonumber\\
 &&+(s^i_{|0}-q^{\ i}_0)y_k+(r_{k0|0}-r_{00|k})b^i,\label{j1}
 \eeq
 where $\lambda=\lambda(x)$ is a scalar function, and $\bar{R}^i_{\ k}$ denotes the Riemann curvature of $\alpha$.
\end{lem}

{\it Proof :} {\bf $\Longrightarrow$ :}
 By the definition of the Weyl curvature
$\bar{W}_{ik}$ of $\alpha$ we have
 \be\label{w52}
 \bar{W}_{ik}=\bar{R}_{ik}-\frac{1}{n-1}\bar{R}ic_{00}a_{ik}+\frac{1}{n-1}\bar{R}ic_{k0}y_i,
 \ee
where $\bar{R}_{ik}:=a_{im}\bar{R}^m_{\ k}$ and $\bar{R}ic_{ik}$
denote the Ricci tensor of $\alpha$. By
$\bar{R}_{ik}=\bar{R}_{ki}$ and (\ref{w52}) we get
 \be\label{w53}
\bar{W}_{ik}-\bar{W}_{ki}=\frac{1}{n-1}\big(\bar{R}ic_{k0}y_i-\bar{R}ic_{i0}y_k\big).
 \ee
  Plugging (\ref{y33}) into
  (\ref{w53}) yields
  \be\label{w54}
 T_ky_i-T_iy_k+(n-1)\big[(s_{k|i}-s_{i|k}+q_{ki}-q_{ik})\alpha^2+(r_{k0|0}-r_{00|k})b_i-(r_{i0|0}-r_{00|i})b_k\big]=0,
  \ee
 where we define
  $$
  T_k:=(n-2)q_{0k}-q_{k0}-(n-1)s_{k|0}-b^lr_{k0|l}-r_{k0}r^l_{\
  l}-\bar{R}ic_{k0}.
  $$
Contracting (\ref{w54}) by $y^kb^i$ we get
 \be\label{w55}
 (\cdots)\alpha^2+\beta[T_0+(n-1)b^l(r_{00|l}-r_{l0|0})]=0.
 \ee
By (\ref{w55})  we obtain
 \be\label{w56}
T_0+(n-1)b^l(r_{00|l}-r_{l0|0})=(n+1)\eta\alpha^2,
 \ee
where $\eta=\eta(x)$ is a scalar function. Then it follows from
the definition of $T_i$ and (\ref{w56}) that
 \be\label{w57}
 \bar{R}ic_{00}=(n-3)q_{00}-(n-1)s_{0|0}-(n+1)\eta\alpha^2+(n-2)b^lr_{00|l}-(n-1)b^lr_{l0|0}-r^l_{\
 l}r_{00}.
 \ee
 By (\ref{w57}) we can get $\bar{R}ic_{k0}$. Plugging
 (\ref{w57}) and $\bar{R}ic_{k0}$ into (\ref{w52}) we get $\bar{W}_{ik}$, and then by (\ref{w52}) and (\ref{y33}) we obtain
 (\ref{j1}), where $\lambda$ is defined by
  \be\label{lam}
  \lambda:=-\frac{(n+1)\eta-s^l_{|l}}{n-1}.
  \ee

  {\bf $\Longleftarrow$ :} Suppose that  (\ref{j1}) holds. Using the first formula in (\ref{j8}) and $b=constant$ we have
 \be
  q_{00}=-s_{0|0}-p_{00}-b^lr_{l0|0}.\label{j2}
 \ee
Contracting  (\ref{j1}) over
  $i,k$ we get $\bar{R}ic_{00}$, and then using (\ref{j2}) we
  obtain (\ref{w57}) with $\eta$ defined by (\ref{lam}). Now plugging (\ref{w57}) and (\ref{j1}) into (\ref{w52}),
  we immediately obtain (\ref{y33}).
    \qed

\

It is clear that no obvious way shows that
 $\bar{R}_{ik}=\bar{R}_{ki}$ in (\ref{j1}). It follows from (\ref{j1}) that
  the symmetric condition
 $\bar{R}_{ik}=\bar{R}_{ki}$ is equivalent to
 \beq
 0&=&\big[b^l(r_{l0|k}+r_{lk|0}-2r_{k0|l})-q_{k0}+q_{0k}-s_{k|0}+s_{0|k}\big]y_i+2(r_{k0|0}-r_{00|k})b_i\nonumber\\
 &&+2(q_{ki}+s_{k|i})\alpha^2-(i/k).\label{y061}
  \eeq

\begin{lem}
(\ref{y32}) and (\ref{j1}) $\Longleftrightarrow $ (\ref{y32}),
(\ref{y48}) and (\ref{y49}).
\end{lem}

{\it Proof :} {\bf $\Longrightarrow$ :} To simplify (\ref{y061}),
we first give two formulas as follows by (\ref{y32}), (\ref{j2})
and (\ref{j1}):
 \beq
 b^l(r_{l0|0}-r_{00|l})&=&\Big[\lambda-\frac{(n-3)s^ls_l-t^l_{\
 l}}{n-1}\Big](\alpha^2-\beta^2)+(t_0-q_0+b^ls_{0|l})\beta,\label{y065}\\
q_{0i}-q_{i0}&=&2\Big[\lambda-\frac{(n-3)s^ls_l-t^l_{\
l}}{n-1}\Big](y_i-\beta b_i)+2(t_0-q_0+b^ls_{0|l})b_i+s_{i|0}-s_{0|i}\nonumber\\
&&-b^l(r_{li|0}+r_{l0|i}-2r_{i0|l}).\label{y066}
 \eeq
To show (\ref{y065}) and (\ref{y066}), by the first formula in
(\ref{j08}) we have
 \be\label{g68}
b^l(r_{l0|0}-r_{00|l})=b^l(s_{l0|0}+b^k\bar{R}_{kl}),\ \ \ \
r_{i0|0}-r_{00|i}=s_{i0|0}+b^k\bar{R}_{ki}.
 \ee
 Contracting (\ref{j1}) by $b_ib^k$, and then using (\ref{j2}),
 the second formula of (\ref{j08}), the first formula of
 (\ref{j008}) and the third formula of (\ref{j9}), we have
 \be\label{g69}
 b^mb^l\bar{R}_{ml}=(\lambda-s_ls^l)\alpha^2-\lambda\beta^2+(b^ls_{0|l}-t_0-q_0)\beta-2s_{0|0}+s_0^2-b^lr_{l0|0}-p_{00}.
 \ee
Similarly, by (\ref{y32}), (\ref{j2}) and  the first formula of
 (\ref{j008}), we have
 \be\label{g70}
 b^ls_{l0|0}=\frac{2s_ls^l+t^l_{\
 l}}{n-1}\alpha^2+\frac{(n-3)s_ls^l-t^l_{\ l}}{n-1}\beta^2+2\beta
 t_0-s_0^2+2s_{0|0}+b^lr_{l0|0}+p_{00}.
 \ee
Then by the first formula in (\ref{g68}), we obtain (\ref{y065})
from (\ref{g69}) and (\ref{g70}). Now by a contraction of
(\ref{j1})
 we get $b^k\bar{R}_{ki}$, and then using the obtained $b^k\bar{R}_{ki}$, (\ref{y32}),
(\ref{y065}) and  the first formula of
 (\ref{j008}), we obtain (\ref{y066}) from the second formula in (\ref{g68}).

Now  contracting (\ref{y061})  by $y^k$ and using (\ref{y065}), we
can write  (\ref{y061})
 as
  \be\label{y067}
 A_i\alpha^2+\beta B_i=0,
  \ee
where $A_i,B_i$ are polynomials in $y$. By (\ref{y067}) we have
$B_i=\sigma_i\alpha^2$, which is expressed as follows
 \be\label{y068}
 r_{00|i}-r_{i0|0}=\Big[\lambda-\frac{(n-3)s^ls_l-t^l_{\
l}}{n-1}\Big]\beta
y_i-(t_0-q_0+b^ls_{0|l})y_i+\frac{\alpha^2}{2(n-1)}\sigma_i.
 \ee
Plugging (\ref{y068}) into (\ref{y067}) yields
 \be\label{y069}
 b^l(r_{li|0}+r_{l0|i}-2r_{i0|l})=\Big[\lambda-\frac{(n-3)s^ls_l-t^l_{\
l}}{n-1}\Big](2y_i-\beta
b_i)+(t_0-q_0+b^ls_{0|l})b_i+\frac{\beta}{2(n-1)}\sigma_i.
 \ee
Now by (\ref{y066}), (\ref{y068}) and (\ref{y069}), we see that
(\ref{y061}) is equivalent to
 \be\label{y070}
 q_{ik}-q_{ki}=s_{k|i}-s_{i|k}+\frac{b_k\sigma_i-b_i\sigma_k}{2(n-1)}.
 \ee

By a contraction on (\ref{y069}), we easily obtain (\ref{j00}) for
the expression of $\sigma_i$  by the second formula of
 (\ref{j08}) and the third formula of (\ref{j9}).
 Now using (\ref{j00}), we can easily obtain (\ref{y48}) by (\ref{j2}) and (\ref{y070}) since
 we can write (\ref{j2}) as $q_{ik}+q_{ki}=\cdots$.

 Finally, by (\ref{j00}),
(\ref{y068}) and  (\ref{y069}) are respectively reduced to
 \beq\label{y071}
r_{00|i}-r_{i0|0}\hspace{-0.5cm}&&=\frac{\alpha^2\sigma_i-\sigma_0y_i}{2(n-1)},\\
b^l(r_{li|0}+r_{l0|i}-2r_{i0|l})\hspace{-0.5cm}&&=2\Big[\lambda-\frac{(n-3)s^ls_l-t^l_{\
l}}{n-1}\Big]y_i+\frac{b_i\sigma_0+\beta\sigma_i}{2(n-1)}.\label{y072}
 \eeq
Now under the formulas (\ref{y48}),  (\ref{y071}) and
(\ref{y072}),
 we can easily show that  (\ref{j1}) is equivalent to (\ref{y49}) with $B_{ik}$ defined by
(\ref{j01}), where we have used (by (\ref{j00}))
 \be\label{g77}
b^l\sigma_l=-2t^l_{\ l}-2(n-1)\lambda+2(n-3)s_ls^l.
 \ee

{\bf $\Longleftarrow$ :} To verify (\ref{j1}), by the last
argument above, we only need to verify (\ref{y071}) and
(\ref{y072}), and then we get (\ref{j1}) following from
(\ref{y49}).

Contracting (\ref{y49}) by $b_ib^k$ and using (\ref{j00}),
(\ref{j2}), (\ref{g77}),  and the second formula of (\ref{j08}),
the first formula of
 (\ref{j008}) and the third and fourth formulas of (\ref{j9}), we  obtain
 \beq\label{g0085}
 b^mb^l\bar{R}_{ml}\hspace{-0.5cm}&&=\Big[2\lambda-\frac{2(n-2)s_ls^l-t^l_{\
 l}}{n-1}\Big](\alpha^2-\beta^2)-s_ls^l\beta^2+2(b^ms_{0|m}-q_0)\beta\nonumber\\
 &&-2s_{0|0}+s_0^2-p_{00}-b^l(2r_{l0|0}-r_{00|l}).
 \eeq
  Similarly, by (\ref{y32}), (\ref{j2}) and the first formula of (\ref{j008}), we have (\ref{g70}). Then by the first formula in (\ref{g68}), we also obtain
 (\ref{y065}) by (\ref{g0085}) and (\ref{g70}). Next we prove (\ref{y071}). First , by (\ref{y48})
 and (\ref{j00}) we have
  \be\label{g79}
 q_{0k}=\frac{\sigma_0b_k-\beta
 \sigma_k}{4(n-1)}-\frac{1}{2}b^l(r_{l0|k}+r_{lk|0})-p_{k0}-s_{0|k}.
  \ee
Now by (\ref{y32}) and (\ref{y49}), we can get  $s_{k0|0}$ and
$b^l\bar{R}_{lk}$ respectively. Then we can obtain (\ref{y071})
from the second formula in (\ref{g68}), by using  (\ref{j00}),
(\ref{y065}), (\ref{g77}), (\ref{g79}), the second formula of
(\ref{j08}), the first formula of (\ref{j008}) and (\ref{j10}).
For (\ref{y072}), it follows from (\ref{y071}) by (\ref{g77}).
\qed

  \

Conversely, let (\ref{y32}), (\ref{y49}), (\ref{y31}) and
(\ref{y48}) be satisfied. Then  $F$ is of scalar flag curvature,
since it is easy to see from the above proof   that the Weyl
curvature of $F$ vanishes if (\ref{y32}), (\ref{y49}), (\ref{y31})
and (\ref{y48})  are satisfied. This completes the proof of
Proposition \ref{prop45}.    \qed

\

{\it Proof of Theorem \ref{th41} :} {\bf $\Longrightarrow$ :} Let
$F$ be of scalar flag curvature. Then we have (\ref{y49}), and
(\ref{y31})--(\ref{y48}) by Proposition \ref{prop45}. It is
obvious that (\ref{y032}) follows from (\ref{y32}) and
(\ref{y48}).

 {\bf $\Longleftarrow$ :} By
Proposition \ref{prop45}, we only need to show that (\ref{y31})
and (\ref{y48}) automatically hold, provided that (\ref{y032}) and
(\ref{y49}) hold. In fact, we can show that (\ref{y032}) directly
implies (\ref{y31}) and (\ref{y48}). By (\ref{j8}) in Lemma
\ref{lem21} and $r_i=-s_i$, we have
 \beq
 &&t_{ik}=-s_{i|k}-s_{k|i}-p_{ik}+\frac{1}{2}b^m(s_{mi|k}+s_{mk|i}-r_{mi|k}-r_{mk|i}),\label{j86}\\
&&q_{ik}=-s_{i|k}-p_{ik}+\frac{1}{2}b^m(s_{mi|k}-s_{mk|i}-r_{mi|k}-r_{mk|i}).\label{j87}
 \eeq
A direct computation from (\ref{y032}) gives
 \beq\label{j88}
 b^ms_{mi|k}\hspace{-0.5cm}&&=\frac{t^l_{\
 l}+2s_ls^l}{n-1}a_{ik}+b_i\Big\{\frac{(n-3)s_ls^l-t^l_{\
 l}}{n-1}b_k-\frac{1}{2}b^mb^v(r_{mk|v}+r_{mv|k})+s^mr_{mk}\nonumber\\
 &&-b^m(s_{k|m}+s_{m|k})\Big\}+b_k\big[t_i+\frac{1}{2}b^mb^v(r_{mv,i}-r_{mi,v})\big]+\frac{1}{2}b^m(r_{mi|k}+r_{mk|i})\nonumber\\
 &&+p_{ik}+s_{i|k}+s_{k|i}-s_is_k.
 \eeq
 Now plugging (\ref{j88}) into (\ref{j86}) and (\ref{j87})
 respectively and using the first formula of (\ref{j008}) and (\ref{j10}), we obtain (\ref{y31}) and
 (\ref{y48}) respectively.

 For the proof of (\ref{y050}), we first get
$\bar{R}ic_{00}$ by (\ref{y49})  and $s^l_{\ 0|l}$ by
(\ref{y032}), and then plugging them into $K=Ric/((n-1)F^2)$
yields (\ref{y050}), where $Ric$ is given by (\ref{ric}) with
$b=1$. \qed

\subsection{Proofs of Theorem \ref{prop40} and Corollary
\ref{prop41}}

{\it Proof of Theorem \ref{prop40} :} \
 It is shown in \cite{Y3} that a Kropina metric
$F=\alpha^2/\beta$ with $||\beta||_{\alpha}=1$ is a Douglas metric
if and only if (\ref{ycw1}) holds. Therefore, by Theorem
\ref{th41}, we only need to use (\ref{ycw1}) to show that
(\ref{y032}) holds. By (\ref{ycw1}) and definitions, we easily get
 \beq\label{ycw79}
 t_{ij}=-s^ls_lb_ib_j-s_is_j, \ \  t_i=-s^ls_lb_i,\ \   t^l_{\
 l}=-2s^ls_l,\ \   q_{ik}=-s_is_k-s^mr_{mi}b_k,\ \
 q_i=s_ls^lb_i.
 \eeq
 Now for the left hand side of (\ref{y032}), we have
  $$
s_{ij|k}\stackrel{(\ref{ycw1})}{=}(r_{ik}+s_{ik})s_j+b_is_{j|k}-(r_{jk}+s_{jk})s_i-b_js_{i|k}
\stackrel{(\ref{ycw1})}{=}s_j(r_{ik}+b_is_k)+b_is_{j|k}-(i/j),
  $$
and for the right hand side of (\ref{y032}), we also obtain the
same result as above by using $t_i,t^l_{\ l}$ and
$q^*_{ik}=q_{ik}$ in (\ref{ycw79}). Thus we have verified
(\ref{y032}). For (\ref{j000}), it directly follows from using the
second  formula in (\ref{j08}) and then plugging $t^l_{\ l},
q_i,t_i$ of (\ref{ycw79}) into (\ref{j00}). \qed

\

\

\noindent{\it Proof of Corollary \ref{prop41} :} \  Since $\beta$
is closed by (\ref{y53}), we see  that (\ref{ycw1}) automatically
holds. Plug (\ref{y53}) into (\ref{j000}) and (\ref{j01}) we get
 \be\label{ycw83}
 \sigma_i=-2(n-1)\lambda b_i,\ \ \ B_{ik}=-\lambda b_ib_k.
 \ee
Now plugging (\ref{y53}) and (\ref{ycw83}) into (\ref{y49}) we
obtain
 \be\label{ycw84}
  \bar{R}^i_{\ k}=-\epsilon^2(\alpha^2\delta^i_k-y^iy_k)+(\lambda+\epsilon^2)(\alpha^2
  b^ib_k+\beta^2\delta^i_k-\beta y^ib_k-\beta y_kb^i).
   \ee
   By (\ref{y53}) and (\ref{ycw83}), it follows from (\ref{y071})
   that
    \be\label{ycw85}
 \lambda+u+\epsilon^2=0.
    \ee
    Then (\ref{ycw84}) and (\ref{ycw85}) imply (\ref{y54}), and we
    get (\ref{y056}) from (\ref{y050}), (\ref{y53}) and
    (\ref{ycw85}).   \qed

\section{Kropina metrics of constant flag curvature}\label{sec5}

It has been solved for the local structure of Kropina metrics of
constant flag curvature (cf. \cite{SYa1} \cite{YO}). In this
section, we will use Theorem \ref{th41} to investigate it.

\begin{cor}\label{cor61}
Let $F=\alpha^2/\beta$ be an $n$-dimensional Kropina
 metric with $||\beta||_{\alpha}=1$.
 Then $F$ is of constant flag curvature if and only if $\alpha$ is
 of constant sectional curvature $\mu$ and $\beta$ satisfies $r_{00}=0$.
 In this case, we have $\mu\ge 0$, and $F$ is flat-parallel ($\alpha$ is flat and $\beta$ is parallel),
 or up to a scaling on $F$, $\alpha$ and $\beta$ can be locally written as
\be\label{w054}
 \alpha=\frac{\sqrt{(1+|x|^2)|y|^2-\langle
 x,y\rangle^2}}{1+|x|^2}, \ \ \ \
\beta=\frac{\langle Ux+e,y\rangle}{1+|x|^2},
 \ee
 where $U=(u^i_j)$ is a skew-symmetric matrix, $e=(e^i)$ is a constant vector
 satisfying
 \be\label{w055}
 |e|=1,\ \ Ue=0, \ \ \delta^{ij}-e^ie^j=\delta^{kl}u^i_ku^j_l.
 \ee
\end{cor}

{\it Proof :} For $n=2$, it has been proved in \cite{SYa1} that
$F$ is flat-parallel. Now assume that $F$ is of constant flag
curvature $K$. Then it follows from Theorem \ref{th41} that its
flag curvature $K$ is given by (\ref{y050}). First we can write
(\ref{y050}) as
 \be\label{w62}
(\cdots)\alpha^2+12(n-1)\beta^4r_{00}^2=0,
 \ee
which implies $r_{00}=c\alpha^2$ for some scalar function
$c=c(x)$. Since $||\beta||_{\alpha}=1$, we have $r_i+s_i=0$. Then
it is easily shown that $c=0$ and hence $r_{00}=0$. Now plug
$r_{ij}=0,r_{ij|k}=0,q_{ij}=0,s_i=0$ into (\ref{w62}) we have
 \be\label{w63}
 (4K-4nK-t^l_{\ l})\alpha^2+4(n\lambda-\lambda+t^l_{\
 l})\beta^2=0.
 \ee
By (\ref{w63}) we easily get
 \be\label{w64}
 K=-\frac{t^l_{\ l}}{4(n-1)}=\frac{\lambda}{4}\ge 0, \ \ ({\rm since}\  t^l_{\ l}\le 0).
 \ee
 Hence we have
  \be\label{g90}
  t^l_{\ l}=-(n-1)\lambda \ (=constant).
  \ee
By $r_{ij}=0,s_i=0$, (\ref{w64}) and (\ref{g90}), it follows from
(\ref{y49}) that
 $$\bar{R}^i_{\ k}=\lambda (\alpha^2\delta^i_k-y^iy_k),$$
which shows that $\alpha$ is of  constant sectional curvature
$\lambda\ge 0$. If $\lambda=0$, then it is easy to show that $F$
is flat-parallel since by (\ref{w64}), we have $t^l_{\ l}=0$ (this
implies that $\beta$ is closed and then parallel by $r_{00}=0$).
If $\lambda>0$, since $r_{00}=0$ and $\alpha$ is of constant
sectional curvature, by solving Killing fields on a Riemannian
space of constant sectional curvature, it follows that, up to a
scaling on $F$, $\alpha$ and $\beta$ can be locally given by
(\ref{w054}) with $U,e$ satisfying (\ref{w055}) (cf. \cite{SYa1}).

Conversely, assume $r_{00}=0$ and $\alpha$ is of  constant
sectional curvature $\mu$ with $||\beta||_{\alpha}=1$. First by
assumption we have
 \be\label{j94}
 r_{ij}=0,\ \ \ \ r_i=0,\ \ \ s_i=0,\ \ \ q_{ij}=r_{im}s^m_{\ j}=0,\ \  \
 t_i=s_ms^m_{\ i}=0.
 \ee
Then by (\ref{j94}), it follows from the first formulas of
(\ref{j08}) and  the second formula of (\ref{j8}) that
 \be\label{j95}
 s_{ij|k}=-b_l\bar{R}^{\ l}_{k\
 ij}=-\mu(b_ia_{jk}-b_ja_{ik}),\ \ \ \ \ \
 t_{ij}=b^ls_{li|j}=-\mu(a_{ij}-b_ib_j),
 \ee
It is clear that the second formula implies $t^l_{\ l}=-(n-1)\mu$.
Now we use Theorem \ref{th41} to verify that $F$ is of constant
flag curvature, namely, we show that (\ref{y032}) and (\ref{y49})
hold for some scalar function $\lambda=\lambda(x)$ and $K$ in
(\ref{y050}) is a constant.  Now define $\lambda:=-t^l_{\
l}/(n-1)=\mu$, and then (\ref{y49}) naturally holds since
$B_{ij}=0$. Finally we verify  (\ref{y032}). By (\ref{j94}), the
first formula of (\ref{j95}) and $t^l_{\ l}=-(n-1)\mu$, we see
that (\ref{y032}) also holds automatically. Therefore, $F$ is of
scalar flag curvature by Theorem \ref{th41}, and its flag
curvature is given by (\ref{y050}). Now by (\ref{y050}) we have
 \be\label{w65}
K=\lambda s^2+\frac{4s^2-1}{4(n-1)}t^l_{\ l}.
 \ee
Since $t^l_{\ l}=-\mu(n-1)$ as shown above, we have
$K=\mu/4=constant$ by (\ref{w65}).  \qed

\section{Construction by warped product method}\label{sec6}
In this section, we will use Corollary \ref{prop41} to show a
family of examples of projectively flat Kropina metrics with
$\alpha$ in warped product.

Let $M=\mathcal{R}\times \widetilde{M}$ be a product manifold,
where $\widetilde{M}$ is an $(n-1)$-dimensional manifold. Let
$\{x^A\}_{A=2}^n$ be a local coordinate
 system on $\widetilde{M}$. A Riemann
metric $\alpha$ of warped product type is defined as
 \be\label{y55}
 \alpha^2=(y^1)^2+h^2(x^1)\widetilde{\alpha}^2,
 \ee
 where $\widetilde{\alpha}^2=\widetilde{a}_{AC}y^Ay^C$ is a Riemann metric on
 $\widetilde{M}$.  The Riemann curvature tensors $\bar{R}$ of
 $\alpha$ and $\widetilde{R}$ of $\widetilde{\alpha}$ in (\ref{y55}) are related by
 \beq
 \bar{R}^1_{\ k}&=&\frac{h''}{h}(y^1y_k-\alpha^2\delta^1_k),\label{y56}\\
\bar{R}^A_{\ C}&=&\widetilde{R}^A_{\
C}-(h')^2(\widetilde{\alpha}^2\delta^A_C-y^A\widetilde{y}_C)-\frac{h''}{h}(y^1)^2\delta^A_C,\label{y57}
 \eeq
where $y_k:=a_{kl}y^l,\widetilde{y}_C:=\widetilde{a}_{CA}y^A$.
Define $\eta=\eta(x^1):=\int h(x^1)dx^1$, and then a direct
computation shows that
 $$\eta_{i|j}=\eta''\alpha^2,\ \ (\eta_i:=\eta_{x^i}),$$
where the covariant derivative is taken with respect to $\alpha$.
The converse is proved in the following.

\begin{lem}(\cite{IT} \cite{Ta})\label{lem51}
 Let $\alpha$ be a Riemann metric on $M$. Suppose there are two
 functions $\eta$ and $\xi$ on $M$ with $d\eta\ne 0$ such that
  $$\eta_{i|j}=\xi a_{ij},\ \ (\eta_i:=\eta_{x^i}).$$
  Then $\alpha$ is a warped product on $M=\mathcal{R}\times \widetilde{M}$, namely,
  locally $\eta$
  depends only on the the parameter $x^1$ of $\mathcal{R}$,
  $\xi=\eta''(x^1)$ and $\alpha$ can be expressed as
$$
 \alpha^2=(y^1)^2+(\eta'(x^1))^2\widetilde{\alpha}^2.
$$
\end{lem}

Now we show  a construction of examples of Kropina metrics of
scalar flag curvature.

\begin{prop}\label{prop51}
 Let $F=\alpha^2/\beta$ be an $n(\ge 2)$-dimensional Kropina metric
on a product manifold $M=\mathcal{R}\times \widetilde{M}$, where
 \be\label{y58}
 \alpha^2:=(y^1)^2+h^2(x^1)\widetilde{\alpha}^2,\ \ \ \beta:=y^1,
 \ee
 where $h\ne 0$ is a smooth function on $\mathcal{R}$ and $\widetilde{\alpha}$ is an $(n-1)$-dimensional Riemann metric on
 $\widetilde{M}$. Then $F$ is locally
 projectively flat if and only
 if $\widetilde{\alpha}$ is locally flat. In this case, the  flag curvature
$K$ is given by
 \be\label{y059}
K=-\big(\frac{\beta}{\alpha}\big)^6\Big\{\frac{h''}{h}+3(h')^2\big(\frac{\widetilde{\alpha}}{\alpha}\big)^2\Big\}.
 \ee
\end{prop}

{\it Proof :} For $n=2$, we can directly verify that
$F=\alpha^2/\beta$ defined by (\ref{y58}) is projectively flat (we
may put $\widetilde{\alpha}=c(x^2)y^2$). We consider $n\ge 3$. For
the $\alpha$ and $\beta$ defined by (\ref{y58}), a direct
computation shows that $||\beta||_{\alpha}=1$ and (\ref{y53})
holds with
 \be\label{ycw98}
 \epsilon=\frac{h'}{h}, \ \ u=\big(\frac{h'}{h}\big)'.
 \ee
So $F$ is locally projectively flat if and only if (\ref{y54})
holds by Corollary \ref{prop41}.

It can be easily verified that (\ref{y54}) is equivalent to
 \be
 \bar{R}^1_{\ k}=\big[-(u+\epsilon^2)(y^1)^2-\epsilon^2 h^2\widetilde{\alpha}^2\big]\delta^1_k
 +(u+\epsilon^2)
 y^1y_k-uh^2\widetilde{\alpha}^2b_k,\label{y59}
 \ee
 and
 \be
 \bar{R}^A_{\ C}=\big[-(u+\epsilon^2)(y^1)^2-\epsilon^2 h^2\widetilde{\alpha}^2\big]\delta^A_C
 +\epsilon^2 h^2y^A\widetilde{y}_C,\label{y60}
 \ee
where $\widetilde{y}_C:=\widetilde{a}_{CA}y^A$. By (\ref{y56}), we
see that (\ref{y59}) is equivalent to
 \be\label{y61}
 u+\epsilon^2=\frac{h''}{h},
 \ee
which automatically holds by (\ref{ycw98}). By the first equation
in (\ref{ycw98}), it follows from (\ref{y57}) that (\ref{y60}) is
equivalent to
 \be\label{y62}
\widetilde{R}^A_{\
C}=\big[-\epsilon^2h^2+(h')^2\big](\widetilde{\alpha}^2\delta^A_C-y^A\widetilde{y}_C)=0.
 \ee

Now suppose $F$ is locally projectively flat. Then we have
(\ref{y62}), namely, $\widetilde{\alpha}$ is locally flat.
Conversely, if $\widetilde{\alpha}$ is locally flat, then by the
above proof, we can easily get (\ref{y54}).

Finally, by (\ref{y056}), we obtain the  flag curvature $K$ given
by (\ref{y059}).          \qed

\

By Proposition \ref{prop51}, $F=\alpha^2/\beta$ in dimension $n\ge
2$ is locally projectively flat, where $\alpha$ and $\beta$ are
defined by (\ref{y58}) with $h\ne 0$ being arbitrary and
$\widetilde{\alpha}$ being locally flat.

\begin{prop}
 Let $F=\alpha^2/\beta$ be an $n(\ge 2)$-dimensional Kropina metric,
 where $\alpha$ and $\beta$ satisfy (\ref{y53}) with $||\beta||_{\alpha}=1$, $d\epsilon\ne 0$
 and $u=f(\epsilon)\ne 0$  for some function $f$. Then $F$ is locally projectively flat
  if and only if $\alpha$ and $\beta$ can be locally written as
 \be\label{y63}
 \alpha^2=(y^1)^2+h^2(x^1)\widetilde{\alpha}^2,\ \ \ \beta=y^1,
 \ee
 where $\widetilde{\alpha}$ is a locally flat Riemann metric. Further,
 $h$  can be actually determined by $f$.
\end{prop}

{\it Proof :} We firstly show (\ref{y63}) by (\ref{y53}). Define
 \be\label{y64}
 \varphi:=\int \frac{1}{f(\epsilon)}\ e^{\int
 \frac{\epsilon}{f(\epsilon)}d\epsilon}d\epsilon.
 \ee
Then by (\ref{y53}) with
 $u=f(\epsilon)\ne 0$, we can easily verify that
  \be\label{y65}
 \varphi_{i|j}=\epsilon \
 e^{\int\frac{\epsilon}{f(\epsilon)}d\epsilon}a_{ij}, \ \
 (\epsilon_i:=\epsilon_{x^i}).
  \ee
Obviously we have $d\varphi \ne 0$. Then by (\ref{y65}) and Lemma
\ref{lem51}, $\alpha$ is a warped product which can be locally
written as the first expression in
 (\ref{y63}) with $h(x^1)=\varphi'(x^1)$. By (\ref{y64}), we can
 define
  $$g(\varphi):=\int \frac{1}{f(\epsilon)}d\epsilon.$$
  Further by (\ref{y53}) we have
  \be\label{y66}
   \beta=\frac{\epsilon_i}{f(\epsilon)}dx^i=\frac{d\epsilon}{f(\epsilon)}
   =d\big(\int
   \frac{1}{f(\epsilon)}d\epsilon\big)=d(g(\varphi))=g'(\varphi)\varphi'(x^1)dx^1.
   \ee
   Then by $||\beta||_{\alpha}=1$, $\alpha$ in (\ref{y63}), and
   (\ref{y66}), we must have $g'(\varphi)\varphi'(x^1)=1$ and
   $\beta=y^1$.

 Therefore,  by
Proposition \ref{prop51}, we conclude that  $F$ is locally
projectively flat if and only
 if $\widetilde{\alpha}$ in (\ref{y63}) is locally flat.  \qed

\section{Proof of Theorem \ref{th2}}

Let $F=\alpha+\beta$ be a Randers metric and $(h,\rho)$ be its
navigation data. It is known that
$$\alpha^2=\frac{(1-b^2)h^2+\rho^2}{(1-b^2)^2},\ \ \ \beta=-\frac{\rho}{1-b^2},\ \ \ (b=||\beta||_{\alpha}=||\rho||_h).$$
By assumption there hold $\lim_{b\rightarrow
1^-}h=\widetilde{\alpha}$ and $\lim_{b\rightarrow
1^-}\rho=\widetilde{\beta}$. Therefore we have
 \beqn
 \lim_{b\rightarrow 1^-}F&=&\lim_{b\rightarrow
 1^-}\Big(\sqrt{\frac{(1-b^2)h^2+\rho^2}{(1-b^2)^2}}
 -\frac{\rho}{1-b^2}\Big)
 =\lim_{b\rightarrow
 1^-}\frac{h^2}{\sqrt{(1-b^2)h^2+\rho^2}+\rho}\nonumber\\
 &=&\frac{\widetilde{\alpha}^2}{2\widetilde{\beta}}=\frac{1}{2}\widetilde{F},
 \ \ \ ( let \   \widetilde{\beta}>0 \  by \  \widetilde{F}>0).
  \eeqn
This proves Theorem \ref{th2} (i).

 Since  $\widetilde{F}=\lim_{b\rightarrow
 1^{-}}2F$, we have $\widetilde{G}^i=\lim_{b\rightarrow
 1^{-}}G^i$.
So for the curvatures $W^i_{\ k},D^{\ i}_{h\ jk},W^o$ of $F$ and
corresponding  $\widetilde{W}^i_{\ k},\widetilde{D}^{\ i}_{h\
jk},\widetilde{W}^o$ of $\widetilde{F}$, we obtain
$$\lim_{b\rightarrow 1^-}\widetilde{W}^i_{\
k}=W^i_{\ k}, \ \ \  \lim_{b\rightarrow 1^-}\widetilde{D}^{\
i}_{h\ jk}=D^{\ i}_{h\ jk},\ \ \  \lim_{b\rightarrow
1^-}\widetilde{W}^o=W^o.
$$
Therefore, if $F$ is of scalar flag curvature (resp. locally
projectively flat, or Douglassian), then $\widetilde{F}$ is also
of scalar flag curvature (resp. locally projectively flat, or
Douglassian).

Now assume $F$ is of weakly isotropic flag curvature with the
 flag curvature $K$  in the form
 $$K=\frac{3\theta}{F}+\sigma,$$
 where $\theta$ is a 1-form and $\sigma=\sigma(x)$ is a scalar
 function. Since $\theta$ and $\sigma$ are uniquely determined by $F$, by taking
 the limit $b\rightarrow 1^-$ on both sides of the above, we see that $\widetilde{F}$ is also of weakly isotropic flag
 curvature. Thus by \cite{SYa1}, $F$ is of constant
 flag curvature. This fact can also be proved in another way.
Let $(h,\rho)$ be the navigation data of the Randers metric
$F=\alpha+\beta$. Since $F$ is of weakly isotropic flag curvature,
it is known that $h=\sqrt{h_{ij}y^iy^j}$ is of isotropic sectional
curvature $\mu=\mu(x)$ (a constant in dimension $n\ge 3$) and
$\rho=\rho_iy^i$ satisfies $\rho_{i|j}+\rho_{j|i}=ch_{ij}$ for
some scalar function $c=c(x)$ (\cite{SYi}). Since
$\lim_{b\rightarrow 1^-}h=\widetilde{\alpha}$ and
$\lim_{b\rightarrow 1^-}\rho=\widetilde{\beta}$, by taking
 the limit $b\rightarrow 1^-$, we have
 $\widetilde{r}_{00}=\widetilde{c}\widetilde{\alpha}^2$ from $\rho_{i|j}+\rho_{j|i}=ch_{ij}$, and
 $\widetilde{\alpha}$ is of isotropic sectional curvature $\widetilde{\mu}$,
 where $\widetilde{c}:=\lim_{b\rightarrow 1^-}c$ and  $\widetilde{\mu}:=\lim_{b\rightarrow 1^-}\mu$. We have
 $\widetilde{c}=0$ from
 $\widetilde{r}_{00}=\widetilde{c}\widetilde{\alpha}^2$ and
 $||\widetilde{\beta}||_{\widetilde{\alpha}}=1$. Further, we have
 $\widetilde{\mu}=\mu=constant$ in dimension $n\ge 3$ and
 in particular $\widetilde{\mu}=0$ in dimension $n=2$. So
 $\widetilde{r}_{00}=0$ and $\widetilde{\alpha}$ is of constant sectional
 curvature. Thus by Corollary \ref{cor61},
 $\widetilde{F}=\widetilde{\alpha}^2/\widetilde{\beta}$ is of
 constant flag curvature.    \qed

\begin{rem}
In Theorem \ref{th2}, let $F=\alpha+\beta$ be a Randers metric and
$(h,\rho)$ be its navigation data. Suppose that
$h=\sqrt{h_{ij}y^iy^j}$ and $\rho=\rho_iy^i$ are locally given by
($\rho^i:=h^{ij}\rho_j$)
 \beqn
 h&=&\frac{\sqrt{(1+\mu |x|^2)|y|^2-\mu \langle
 x,y\rangle^2}}{1+\mu |x|^2},\\
 \rho^i&=&-2(\lambda \sqrt{1+\mu |x|^2}+\langle
 d,x\rangle)x^i+\frac{2|x|^2d_i}{1+\sqrt{1+\mu
 |x|^2}}+u^i_kx^k+e^i+\mu \langle e,x\rangle x^i,
  \eeqn
where $\lambda,\mu$ are constants, $U=(u^i_k)$ is a skew-symmetric
matrix and $d,e\in \mathcal{R}^n$ are constant vectors. To take
$b=||\beta||_{\alpha}\rightarrow 1^-$, we only require
$h_{ij}\rho^i\rho^j=1$. A direct computation gives a Kropina
metric $\widetilde{F}=\widetilde{\alpha}^2/\widetilde{\beta}$ in
two cases: (A). $\widetilde{\alpha}=|y|,  \
\widetilde{\beta}=\langle e,y\rangle$; (B). $\widetilde{\alpha}=h$
and $\widetilde{\beta}$ is given by
 $$\widetilde{\beta}=\frac{\langle Ux+e,y\rangle}{1+|x|^2},$$
 where $U$ and $e$
 satisfy
 $$
 |e|=1,\ \ Ue=0, \ \ \mu(\delta^{ij}-e^ie^j)=\delta^{kl}u^i_ku^j_l.
 $$
\end{rem}

\begin{rem}\label{rem72}
For a given Kropina metric $F=\alpha^2/\beta$ with
$||\beta||_{\alpha}=1$, we can construct a family of Randers
metrics $\bar{F}=\bar{\alpha}+\bar{\beta}$ with $F$ as the  limit
of $\bar{F}$ as $\bar{b}=||\bar{\beta}||_{\bar{\alpha}}\rightarrow
1^-$. Define
 $$\bar{\alpha}^2=\frac{(1-\bar{b}^2)\alpha^2+\bar{b}^2\beta^2}{(1-\bar{b}^2)^2},\ \ \
 \bar{\beta}=-\frac{\bar{b}\beta}{1-\bar{b}^2},\ \ \ (|\bar{b}|<1).$$
 Then it can be easily verified that
 $\bar{F}=\bar{\alpha}+\bar{\beta}$ is a Randes metric with
 $\bar{b}=||\bar{\beta}||_{\bar{\alpha}}$ and $F=\lim_{\bar{b}\rightarrow
 1^{-}}2\bar{F}$.
\end{rem}

\vspace{0.5cm}

\noindent Guojun Yang \\
Department of Mathematics \\
Sichuan University \\
Chengdu 610064, P. R. China \\
 yangguojun@scu.edu.cn

\end{document}